\documentclass{article}%
\usepackage{amssymb,amsmath,amsthm}%
\providecommand{\U}[1]{\protect\rule{.1in}{.1in}}
\newcommand\bkE{{\mathbb {E}}}
\newcommand\E{{\mathbb {E}}}
\newcommand\Var{{\mathrm {Var}}}
\newtheorem {Lemma}{Lemma}[section]

\newtheorem {Theorem}{Theorem}[section]
\newtheorem {Proposition}{Proposition}[section]

\newtheorem{Remark}{Remark}[section]

\newcommand\beq{\begin{equation}}
\newcommand\eeq{\end{equation}}

\begin{document}
\begin{center} {\bf \Large Rates in the strong invariance principle for ergodic automorphisms of the torus}\vskip15pt

J\'er\^ome Dedecker $^{a}$, Florence Merlev\`{e}de $^{b}$ and Fran\c{c}oise P\`ene  $^{c}$ \footnote{%
Supported in part by the ANR project PERTURBATIONS}
\end{center}
\vskip10pt
$^a$ Universit\'e Paris Descartes, Sorbonne Paris Cit\'e, Laboratoire MAP5
and CNRS UMR 8145. Email: jerome.dedecker@parisdescartes.fr\\ \\
$^b$ Universit\'e Paris Est, LAMA and CNRS UMR 8050.\\
E-mail: florence.merlevede@univ-mlv.fr\\ \\
$^c$ Universit\'e de Brest,  Laboratoire de Math\'ematiques de Bretagne
Atlantique UMR CNRS 6205.
E-mail: francoise.pene@univ-brest.fr\vskip10pt
{\it Key words}: Invariance principles, strong approximations,  ergodic automorphisms of the torus.\vskip5pt


\begin{center}
{\bf Abstract}\vskip10pt
\end{center}
Let $T$ be an ergodic automorphism of the $d$-dimensional torus ${\mathbb T}^d$.
In the spirit of Le Borgne  \cite{SLB},
we give conditions on the Fourier coefficients of a function $f$ from ${\mathbb T}^d$
to ${\mathbb R}$ under which the partial sums  $f\circ T+ f\circ T^2 +
\cdots + f \circ T^n$ satisfies a strong invariance principle. Next, reinforcing  the condition on the Fourier coefficients
in a natural way, we obtain explicit rates of convergence in the strong invariance principle, up to $n^{1/4} \log n$.

\section{Introduction}
 We endow the $d$-dimensional torus ${\mathbb T}^d={\mathbb R}^d/{\mathbb Z}^d$
with the Lebesgue measure $\bar\lambda$, and we denote by $\mathbb E(\cdot)$ the expectation
with respect to $\bar\lambda$. As usual, the ${\mathbb L}^{p} $ norm of a  $f$ from ${\mathbb T}^d$ to ${\mathbb R}$ is denoted by
$\|f \|_{p}= (  \E ( | f|^p)  )^{1/p}$.

For $d\geq 2$, let $T$ be an ergodic automorphism of
 ${\mathbb T}^d$, and let $f$ be a function from ${\mathbb T}^d$ to ${\mathbb R}$
 such that ${\mathbb E}(f^2)< \infty$ and ${\mathbb E}(f)=0$. In
 \cite{SLB}, Le Borgne has proved that if the Fourier coefficients
 $(c_{\mathbf k})_{\mathbf k\in\mathbb Z^d}$
 of $f$ are such that, for $\theta>2$ and every integer $b>1$,
 \beq \label{Scond}
 \sum_{|\mathbf k|\ge b}|c_{\mathbf k }|^2\le R \log^{-\theta}(b) \, , \quad \text{
 where $|{\bf k}|=\max_{ 1 \leq i \leq d} |k_i|$} \, ,
 \eeq
 then the  partial sums process
 \beq \label{PSP}
  \Big \{ \sum_{i=1}^{[nt]} f \circ T^i, t \in [0,1] \Big \}
 \eeq
properly normalized,
satisfies both the weak and strong invariance principles. More
precisely, Le Borgne has  introduced  in $\cite{SLB}$ an appropriate
$\sigma$-field ${\mathcal F}_0$ such that ${\mathcal F}_0 \subseteq T^{-1}({\mathcal F}_0)$, for which
the quantities $\|{\mathbb E}(f\circ T^k|{\mathcal F}_0)\|_2$ and
$\|f\circ T^{-k}-{\mathbb E}(f\circ T^{-k}|{\mathcal F}_0)\|_2$ can be controlled
for any positive integer $k$.
The weak and strong invariance principles follow then, by applying  Gordin's result
(see \cite{Go}) and
Heyde's result (see \cite{Heyde2}) respectively.

In  Theorem \ref{asip} of this paper, we show that the  weak and strong invariance principles
still hold for functions $f$ satisfying (\ref{Scond}) with $\theta>1$ only,
and we  give a multivariate version of these results.
For the weak invariance principle, this follows from an improvement
of Gordin's criterion, which was already known in the univariate case (see
\cite{DMV}). For the strong invariance principle, this will follow
from a new criterion for stationary sequences, presented  in Theorem
\ref{ASIPgene}
of the appendix.
Note that the condition (\ref{Scond}) with $\theta>1$ is satisfied if, for a positive constant $A$,
\beq \label{NEW}
|c_{\mathbf k}|^2 \leq A \prod_{i=1}^d \frac{1}{(1+|k_i|)
  \log^{1 + \alpha}(2+|k_i|)} \quad \text{for some $\alpha>1$,}
\eeq
 improving on the condition $\alpha>2$
given by Leonov in 1969 (see \cite{Leonov}, Remark 1). Note that  Leonov
has also given
a condition in terms of the modulus of continuity of $f$ in ${\mathbb L}^2$.

The strong invariance principle means that, enlarging ${\mathbb T}^d$
if necessary, there exists a sequence of independent identically distributed (iid)
Gaussian random variables $Z_i$ such that
\beq \label{resautoasip}
\sup_{1\leq k \leq n} \Big|\sum_{i=1}^kf \circ T^i - \sum_{i=1}^k Z_i\Big|  = o \big (  n^{1/2}( \log \log  n)^{1/2}  \big )
\text{ almost surely, as $n\rightarrow \infty$}.
\eeq
 It is also possible to exhibit rates of convergence in
 \eqref{resautoasip}, provided that we reinforce the assumption \eqref{Scond}.
This has been done recently,  thanks to a general result giving  rates of convergence in the strong invariance principle for partial sums of stationary sequences. More precisely, let $p\in]2,4]$ and  $q=p/(p-1)$. We have proved in Theorem 2.1 of \cite{DMP1} that if there exists $R >0$ such that
for every integer $b>1$,
\beq \label{condF1*}\sum_{|\mathbf k|\ge b}|c_{\mathbf k}|^q\le R \log^{-\theta}(b) \ \text{ for some $\theta >\frac {p^2-2}{p(p-1)}$} \, ,\eeq
and
\beq \label{condF2*} \sum_{|\mathbf k|\ge b}|c_{\mathbf k}|^2\le R \log^{-\beta}(b) \ \text{ for some $\beta >\frac {3p-4}{p}$} \, ,\eeq
then the strong approximation \eqref{resautoasip} holds true with
an error of order $o \big ( n^{1/p} ( \log n)^{(t+1)/2}  \big )$, for $t >2/p$.
A condition on the $\ell^q$ norm of $(c_{\mathbf k})_{\mathbf k\in\mathbb Z^d}$ seems
appropriate in this context, since this $\ell^q$-norm dominates  the ${\mathbb L}^p$ norm
of $f$, which is required to be finite to get the rate $o(n^{1/p})$ in the iid situation.

If we assume that the Fourier coefficients of $f$ are such that,
\begin{equation}\label{Leonov}
  |c_{\mathbf k}|^q \leq A \prod_{i=1}^d \frac{1}{(1+|k_i|)
  \log^{1 + \alpha}(2+|k_i|)} \, ,
\end{equation}
then the conditions \eqref{condF1*} and \eqref{condF2*} are both satisfied provided that $\alpha >(p^2-2)/(p^2-p)$. Now, considering (\ref{NEW}), one can wonder if $\alpha >1$
in (\ref{Leonov}) is enough to  get an approximation error  of order $o(n^{1/p} L(n))$ in \eqref{resautoasip}, where $L(n)$ is a slowly varying function. The main result
of this paper, Theorem \ref{automorphismes} below, shows that the answer is positive.

\setcounter{equation}{0}

\section{Invariance principles for ergodic  automorphisms of the torus}
\setcounter{equation}{0}

Let us first recall some probabilistic notations.
 A measurable function $f:{\mathbb T}^d\rightarrow {\mathbb R}^m$
(with coordinates $f_1,...,f_m$)
is said to be centered if every $f_i$ is integrable and centered.
Such a function $f$ is said to be square integrable if every $f_i$ is square integrable.
Now, for every centered and square integrable functions $f,g:{\mathbb T}^d\rightarrow
{\mathbb R}^m$ (with $f=(f_1,...,f_m)$ and
$g=(g_1,...,g_m)$), we define the covariance matrix
$\text{Cov}(f,g)$ of $f$ and $g$ and the variance matrix $\text{Var}(f)$ by
$$\text{Cov}(f,g)=(\mathbb E(f_ig_j))_{i,j=1,...,m}\, , \quad \text{and} \quad \text{Var}(f)= \text{Cov}(f,f)\, .$$

Let us now recall some facts about ergodic automorphisms of $\mathbb T^d$.
A group automorphism $T$ of $\mathbb T^d$ is the quotient map of a linear map
$\tilde T:{\mathbb R}^d\rightarrow{\mathbb R}^d$ given by $\tilde T(x)=S.x$
($.$ being the matrix product), where $S$ is a $d\times d$-matrix
with integer entries and with determinant $\pm1$.
Any automorphism $T$ of $\mathbb T^d$ preserves the Lebesgue measure $\bar\lambda$.
Therefore $(\mathbb T^d,{\mathcal B}(\mathbb T ^d),T,\bar\lambda)$ is
a probability dynamical system (where ${\mathcal B}(\mathbb T ^d)$ stands for the Borel
$\sigma$-algebra of $\mathbb T^d$).

This dynamical system is ergodic if and only if no root of the unity is an eigenvalue of
the matrix $S$ associated to $T$. In this case, we say that $T$ is an ergodic
automorphism of $\mathbb T^d$.

An automorphism $T$ of $\mathbb T^d$ is said to be hyperbolic if the matrix $S$ associated to
$T$ admits no eigenvalue of modulus one. With the preceding characterization of ergodic
automorphisms of $\mathbb T^d$, it is clear that every hyperbolic automorphism of $\mathbb T^d$
is ergodic.
Ergodic automorphisms of $\mathbb T^d$ are partially hyperbolic but not necessarily
hyperbolic (an example of a non-hyperbolic ergodic automorphism of $\mathbb T^d$
can be found in \cite{SLB}).

\medskip

In the next Theorem, we  give  weak and strong invariance principles for the partial
sum process (\ref{PSP}) of
${\mathbb R}^m$-valued functions.

\begin{Theorem}\label{asip}
Let $T$ be an ergodic automorphism of $\mathbb T^d$. For any $j\in \{ 1, \ldots, m \}$, let $f_j:\mathbb T^d\rightarrow\mathbb R$ be a centered function and assume that its  Fourier coefficients $(c_{\mathbf k, j })_{\mathbf k\in\mathbb Z^d}$ satisfy the following condition: there exists a positive constant $R$ such that for every integer $b >1 $,
\beq \label{condFasip}\sum_{|\mathbf k|\ge b}|c_{\mathbf k , j}|^2\le R \log^{-\theta}(b) \
\text{ for some }\theta >1 \, .\eeq Let $f=(f_1, \ldots, f_m):\mathbb T^d\rightarrow\mathbb R^m$.
Then the series
$
 \Sigma = \sum_{k \in {\mathbb Z}}  {\rm Cov} (f , f \circ T^k)
$
converges, and
\begin{equation}\label{defsigma2}
\lim_{n \rightarrow \infty} \frac{1}{n}{\rm{Var}}\Big(\sum_{i=1}^n f \circ T^i\Big)= \Sigma
\, .
\end{equation}
 In addition,
\begin{itemize}\item[1.] The process $ \{ n^{-1/2} \sum_{i=1}^{[nt]} f \circ T^i, t \in [0,1] \}$ converges in $D([0,1], {\mathbb R}^m )$ equipped with the uniform topology to a Wiener process $\{W(t), t\in [0,1]\}$ with variance matrix ${\rm Var} (W(1)) =\Sigma $.
\item[2.] Enlarging ${\mathbb T}^d$ if
necessary, there exists a sequence $(Z_i)_{i
\geq 1}$ of iid   ${\mathbb R}^m$-valued
Gaussian random variables  with zero mean and
variance matrix ${\rm Var}(Z_1)=\Sigma$ such that
$$
\sup_{1\leq k \leq n} \Big|\sum_{i=1}^kf \circ T^i - \sum_{i=1}^k Z_i\Big|  = o \big (  n^{1/2}( \log \log  n)^{1/2}  \big )
\text{ almost surely, as $n\rightarrow \infty$}.
$$
\end{itemize}
\end{Theorem}

When $m=1$, it is also possible to exhibit rates of convergence in \eqref{resautoasip} provided that we reinforce Condition \eqref{condFasip}.

\begin{Theorem}\label{automorphismes}
Let $T$ be an ergodic automorphism of $\mathbb T^d$. Let $p\in]2,4]$ and $q:=p/(p-1)$.
Let $f:\mathbb T^d\rightarrow\mathbb R$ be a centered function with Fourier coefficients
$(c_{\mathbf k})_{\mathbf k\in\mathbb Z^d}$ satisfying the following conditions: there exists a positive constant $R$ such that for every integer $b >1 $,
\beq \label{condF1}\sum_{|\mathbf k|\ge b}|c_{\mathbf k}|^q \le R \log^{-\theta}(b) \
\text{ for some }\theta >1 \, ,\eeq
and
\beq \label{condF2} \sum_{|\mathbf k|\ge b}|c_{\mathbf k}|^2 \le  R \,  b^{-\zeta}\ \ \
\text{for some }\zeta >0 \, .\eeq
Then the series
\beq \label{defsigma2b} \sigma^2= \sum_{k\in\mathbb Z}{\mathbb E}(f.f\circ T^k)
\eeq
converges absolutely and, enlarging ${\mathbb T}^d$ if
necessary, there exists a sequence $(Z_i)_{i
\geq 1}$ of iid Gaussian random variables with zero mean and
variance $\sigma^2$ such that \beq \label{resauto}
\sup_{1\leq k \leq n} \Big|\sum_{i=1}^kf \circ T^i - \sum_{i=1}^k Z_i\Big|  = o \big ( n^{1/p}  \log n  \big )
\text{ almost surely, as $n\rightarrow \infty$}.
\eeq
\end{Theorem}
Observe that if \eqref{Leonov} holds with $\alpha >1$ then  \eqref{condF1} and \eqref{condF2} are both satisfied, so that the strong approximation
(\ref{resauto}) holds. However Theorem 2.1 in \cite{DMP1} and Theorem \ref{automorphismes} above have different ranges of applicability. Indeed,
let  $\gamma >1$, and define
 $c_{\bf k} = \ell^{- \gamma/q }$ if ${\bf k} = (2^{\ell}, 0, \ldots, 0)$,
$c_{\bf k} = -\ell^{- \gamma/q }$ if ${\bf k} = (-2^{\ell}, 0, \ldots, 0)$
for $\ell \in {\mathbb N}$,
and $c_{\bf k}  = 0$ otherwise. Let now $b$ and $r$ be positive integers such that $2^{r-1} < b \leq 2^r$.  Since
$$
\sum_{|{\bf k}|\ge b}|c_{\bf k}|^q =
2\sum_{\ell \geq r} \frac{1}{\ell^\gamma} \, ,
$$ it follows that  $\lambda_1 (\log b)^{1- \gamma} \leq \sum_{|{\bf k}| \ge b}|c_{\bf k}|^q  \leq \lambda_2 (\log b)^{1- \gamma} $ (where $\lambda_1$ and $\lambda_2$ are two positive constants). 
Similarly $\lambda_1 (\log b)^{1- 2 \gamma /q } \leq  \sum_{|{\bf k}| \ge b}|c_{\bf k}|^2 \leq \lambda_2 (\log b)^{1- 2 \gamma /q } $. In this situation, the conditions \eqref{condF1*} and \eqref{condF2*} are both satisfied provided that $ \gamma > 1 + (p^2-2)/(p^2-p)$ whereas condition \eqref{condF2} fails.

\medskip

To prove Theorem \ref{automorphismes}, we shall still use martingale approximations as done in \cite{DMP1}, but with the following modifications:
Condition \eqref{condF2} allows us to consider a non stationary sequence $X_\ell^*=f_\ell\circ T^\ell$,  where the functions $f_\ell$ are defined  through a truncated series of the   Fourier coefficients
of $f$. For the partial sums associated to this  non stationary sequence, the approximation error by a non stationary martingale  can be  suitably handled with the help of  Condition \eqref{condF1}.

\section{Proofs of Theorems \ref{asip} and \ref{automorphismes}} \label{sectionproofauto}
\setcounter{equation}{0}
As in \cite{DMP1}, we consider the filtration as defined in \cite{Lind,SLB} that enables to suitably approximate the partial sums $\sum_{i=1}^nf \circ T^i $ by a martingale. To be more precise, given a finite partition $\mathcal P$ of $\mathbb T^d$, we define
the measurable partition $\mathcal P_0^{\infty}$ by~:
$$\forall \bar x\in\mathbb T^d,\ \
 \mathcal P_0^{\infty}(\bar x):=\bigcap_{k\ge 0}T^k\mathcal P(T^{-k}(\bar x))$$
and, for every integer $n$, the $\sigma$-algebra $\mathcal F_n$ generated by
$$\forall \bar x\in\mathbb T^d,\ \
  \mathcal P_{-n}^{\infty}(\bar x):=\bigcap_{k\ge -n}T^k\mathcal P(T^{-k}(\bar x))
      =T^{-n}(\mathcal P_0^{\infty}(T^n(\bar x)) \, .$$
These definitions coincide with the ones of \cite{SLB} applied to
the ergodic toral automorphism $T^{-1}$.
We obviously have $\mathcal F_n\subseteq\mathcal F_{n+1}=T^{-1}\mathcal F_n$. Note that the sequence $(f \circ T^i)_{i \geq 1}$ is non adapted to $({\mathcal F}_i)_{i \geq 1}$.

In what follows, we use the notation $\E_n (f) = \E ( f |{\mathcal F_n})$.

\subsection{Proof of Theorem \ref{asip}}
According to Theorem \ref{ASIPgene} and Remark \ref{use} given in Appendix,
it suffices to verify that condition \eqref{usecond2} is satisfied.
Therefore, it suffices to verify that for any $j \in \{1, \ldots, m \}$,
\beq \label{p1asip}
\sum_{n\geq 3 }\frac{\log n}{n^{1/2}(\log \log n)^{1/2}}\Vert \E_0 (f_j \circ T^n)\Vert _{2}<\infty \text{ and }\sum_{n\geq 3 }\frac{\log n}{n^{1/2}(\log \log n)^{1/2}}%
\Vert f_j-\E_n (f_j) \Vert _{2}<\infty \, .
\eeq
But, according to the proof of Propositions 4.2 and 4.3 of \cite{DMP1} (see also \cite{SLB}), for any $f_j$ satisfying \eqref{condFasip},
\begin{equation*}
\Vert {\mathbb E}_{-n}(f_j)\Vert_2+\Vert f-{\mathbb E}_{n}(f_j)\Vert_2
    \ll n^{-\theta/2} \, .
\end{equation*}
Since $\theta >1$, \eqref{p1asip} is satisfied.  $\square$

\subsection{Proof of Theorem \ref{automorphismes}}

 Let $f:\mathbb T^d\rightarrow \mathbb R$ be a
centered function with  Fourier coefficients
$(c_{\mathbf k})_{\mathbf k\in\mathbb Z^d}$.
For every nonnegative integer $m$, we write
\begin{equation}f_m:=\sum_{|k|\le m}c_k e^{2i\pi\langle k,\cdot\rangle} \, .
\end{equation}
Notice that if $f$ satisfies \eqref{condF2}, then
\begin{equation} \label{majdiff2}
\Vert f-f_m\Vert_2\le Rm^{-\zeta/ 2} \, ,
\end{equation}
and if $f$ satisfies \eqref{condF1}, then
\begin{equation} \label{majdiffp}
\Vert f-f_m\Vert_p\le R ( \log(m) )^{- \theta (p-1)/p} \, .
\end{equation}

According to the proofs of Propositions 4.2, 4.3 and 4.4 of \cite{DMP1},
there exist $c\ge 1$ and
$\gamma,\lambda\in(0,1)$ such that,  setting $b(n):=[\gamma^{-n}]$, we have
\begin{equation}\label{auto1}
\sup_{m\le  b(n)}(\Vert {\mathbb E}_{-n}(f_m)\Vert_p+\Vert f_m-{\mathbb E}_{n}(f_m)\Vert_p)
\ll \lambda^n
\end{equation}
(according to (4.50), (4.51) and (4.53) of \cite{DMP1}),
and
\begin{equation}\label{auto2}
\sup_{N\ge cn}\sup_{m\le b(n)}\sup_{\ell\in\{0,...,n\}}\Vert {\mathbb E}_{-N}(f_mf_m\circ T^\ell)
      -{\mathbb E}(f_mf_m\circ T^\ell)\Vert_{p/2}\ll\lambda^n
\end{equation}
(according to (4.61) and (4.62) of \cite{DMP1}).
Moreover, according to the proof of Propositions 4.2 and 4.3 of \cite{DMP1},
we have, for any $f$ satisfying \eqref{condF1},
\begin{equation}\label{auto3}
\sup_{m\ge 1}(\Vert {\mathbb E}_{-n}(f_m)\Vert_p+\Vert f_m-{\mathbb E}_{n}(f_m)\Vert_p)
    \ll n^{-\theta(p-1)/p} \, ,
\end{equation}
and
\begin{equation}\label{auto4}
\Vert {\mathbb E}_{-n}(f)\Vert_p+\Vert f-{\mathbb E}_{n}(f)\Vert_p
    \ll n^{-\theta(p-1)/p} \, .
\end{equation}
In addition, according to the proof of Proposition 4.4 of \cite{DMP1}, there exists a positive
integer $c$, such that for any $f$ satisfying \eqref{condF1} and \eqref{condF2},
\begin{equation}\label{auto5}
\max_{1 \leq k \leq n}\Vert \E_{-n c }(S^2_k(f)) - \E(S_k^2(f)) \Vert_{p/2} \ll n^{2 -2 \theta (p-1)/p} \, .
\end{equation}

For any $f$ satisfying \eqref{condF2}, using the arguments developed in the proofs of Propositions 4.2 and 4.3 of \cite{DMP1},
we infer that there exists $\beta \in (0,1)$ such that
\begin{equation}\label{auto6}
\sup_{m\ge 1}( \Vert {\mathbb E}_{-n}(f_m)\Vert_2+\Vert f_m-{\mathbb E}_{n}(f_m)\Vert_2 )
    \ll \beta^{n} \, ,
\end{equation}
and
\begin{equation}\label{auto6bis}
 \Vert {\mathbb E}_{-n}(f)\Vert_2+\Vert f-{\mathbb E}_{n}(f)\Vert_2
    \ll \beta^{n} \, .
\end{equation}
Let us write $P_\ell(\cdot)={\mathbb E}_\ell(\cdot)-{\mathbb E}_{\ell-1}(\cdot)$.
Now, let $\alpha$ be a positive real such that $\alpha \zeta \geq 3 -2/p   $. We then define
$$d_1^*:=\sum_{k\in\mathbb Z}P_1\left(f_{1}\circ T^k\right) \ , \ X_1^*:=f_1\circ T \, ,$$
and, for every $j\ge 0$ and every $\ell\in\{2^{j}+1,...,2^{j+1}\}$,
$$d^*_\ell:=\sum_{k\in\mathbb Z}P_\ell\left(f_{[2^{\alpha j}]}\circ T^k\right) \  , \
        X_\ell^*:=f_{[2^{\alpha j}]}\circ T^\ell.$$
For every positive integer $n$, we define
$$M^*_n(f):=\sum_{\ell=1}^n d^*_\ell\ \ \mbox{and}\ \ S_n^*(f):=\sum_{\ell=1}^n X_{\ell}^* \, .$$
The conclusion of Theorem \ref{automorphismes} comes from the three following lemmas.
\begin{Lemma}\label{Lem1}
We have
$|S_n(f)-S_n^*(f)|=o(n^{1/p}(\log n))$ almost surely.
\end{Lemma}
\begin{Lemma}\label{Lem2}
We have
$|S_n^*(f)-M_n^*(f)|=o(n^{1/p}(\log n))$ almost surely.
\end{Lemma}
\begin{Lemma}\label{Lem3}
The conclusion of Theorem  \ref{automorphismes} holds with $M_n^*(f)$ replacing $S_n(f)$.
\end{Lemma}
\noindent{\bf Proof of Lemma \ref{Lem1}.}
For any nonnegative integer $j$, let
\begin{eqnarray*} \label{0dec}
D_j := \sup_{ 1 \leq k \leq 2^{j}} \Big|
\sum_{\ell= 2^{j}+1}^{k+2^{j}}(X_\ell-X_\ell^*)\Big| \, .
\end{eqnarray*}
Let $N \in {\mathbb N}^*$ and let $k \in ]1, 2^N]$. We first notice that
$D_j \geq |  \sum_{\ell= 2^{j}+1}^{2^{j+1}}(X_\ell-X_\ell^*)|$, so if $K$ is the  integer
such that $2^{K-1} < k \leq  2^{K}$, then
$$\big  | S_k-S_k^* \big | \leq |X_1-X_1^*|+\sum_{j=0}^{K-1} D_j\, .$$
Consequently, since $K \leq N$,
\begin{equation} \label{1dec} \max_{1\le k\le 2^N}|S_k-S_k^*|\le |X_1-X_1^*|+\sum_{j=0}^{N-1}
   D_j\, .
    \end{equation}
Therefore, by standard arguments, Lemma \ref{Lem1} will follow if we can prove that
$
D_j= o\left(j \, 2^{j/p} \right)
$ almost surely.
This will hold true as soon as
\begin{equation} \label{lma1p1}
\sum_{j\ge 1}\frac{ \Vert D_j \Vert_q^q }{2^{jq /p}\, j^{q}}<\infty \, \text{ for some $q \in [1,p]$}\, .
\end{equation}
We shall verify (\ref{lma1p1}) for $q=2$. Notice that
\[
 \Vert D_j  \Vert_2 = \Big  \Vert  \sup_{1\le k\le 2^{j}}
    \Big |\sum_{\ell= 2^{j}+1}^{k+2^{j}}(f - f_{[2^{\alpha j}]})\circ T^{\ell}\Big |\Big \Vert_2
    \leq 2^{j} \Vert f - f_{[2^{\alpha j}]} \Vert_2 \, .
\]
Hence, by using \eqref{condF2},
$
\Vert D_j \Vert^2_2   \ll 2^{2j} 2^{-\zeta \alpha j}
$,
which together with the fact that $\alpha \zeta \geq 2 -2/p$  implies \eqref{lma1p1} with $q=2$, and then Lemma \ref{Lem1}. $\square$

\bigskip

\noindent{\bf Proof of Lemma \ref{Lem2}.}  Without loss of generality, we assume that $\theta<(p^2-2)/(p(p-1))$.
Following the beginning of the proof of Lemma
\ref{Lem1}, Lemma \ref{Lem2} will be proven if   \eqref{lma1p1} holds
with $D_j$ defined by
\begin{eqnarray} \label{0declma2}
D_j = \sup_{ 1 \leq k \leq 2^{j}} \Big|  \sum_{\ell= 2^{j}+1}^{k+2^{j}}(X_\ell^* -d_{\ell}^*)\Big| \, .
\end{eqnarray}
With this aim, setting, for every $k\in\{1,...,2^j\}$,
$$
R_{j,k}=  \sum_{\ell= 1}^{k} \Big ( f_{[2^{\alpha j}]}\circ T^\ell  -
   \sum_{m\in\mathbb Z}P_\ell\left(f_{[2^{\alpha j}]}\circ T^m\right) \Big )
   =\sum_{\ell=2^j+1}^{k+2^j}(X_\ell^*-d_\ell^*)\circ T^{-2^j}\, ,
$$
 we first observe that
 \begin{eqnarray}  \label{1declma2} \Vert D_j\Vert_p =
 \Big \Vert \sup_{ 1 \leq k \leq 2^{j}} | R_{j,k}| \Big \Vert_p
 \ll 2^{j/p} \sum_{k=0}^{j} 2^{-k/p} \Vert R_{j,2^k} \Vert_p \ ,
\end{eqnarray}
(where for the  inequality we have used inequality (6) in \cite{W}).
Now, according to the proof of Proposition 5.1 in \cite{DMP1} with
$X_{\ell}= f_{[2^{\alpha j}]} \circ T^{\ell}$ and using again stationarity, we get that for any integer $k \geq 0$ and any integer $N \geq 2^k$,
\begin{multline} \label{2declma2}
\max_{1 \leq m \leq 2^k}\Vert R_{j,m} \Vert_p  \ll  \sum_{\ell=1}^N  \Vert \E_{-\ell }
(f_{[2^{\alpha j}]} ) \Vert_p  +     \sum_{\ell=0}^{N-1}  \Vert f_{[2^{\alpha j}]}- {\mathbb E}_{\ell} (f_{[2^{\alpha j}]} )  \Vert_p \\
 + \Big (  \sum_{m=1}^{2^k}
\big \Vert \sum_{\ell \geq m +N }   P_{-\ell}(f_{[2^{\alpha j}]} )\circ T^\ell \big \Vert^2_p  \Big )^{1/2}
+ \Big (  \sum_{m=1}^{2^k}
  \big \Vert  \sum_{\ell \geq m +N }  P_{\ell}(f_{[2^{\alpha j}]} )\circ T^{-\ell}
   \big \Vert^2_p  \Big )^{1/2} \, .
\end{multline}
Let us first consider the case where $[2^{\alpha j}] \leq b(2^k)$.  Starting from (\ref{2declma2}) with $N=2^k$ and using the fact that $ \Vert P_{-\ell}(f_{[2^{\alpha j}]} ) \Vert_p \leq 2 \Vert \E_{-\ell}(f_{[2^{\alpha j}]}  ) \Vert_p$ and that $\Vert P_{\ell}(f_{[2^{\alpha j}]}  ) \Vert_p  \leq \Vert f_{[2^{\alpha j}]}  - \E_{\ell-1}(f_{[2^{\alpha j}]}  ) \Vert_p$, we get that
\begin{multline*}
\max_{1 \leq m \leq 2^k}\Vert R_{j,m} \Vert_p \ll  \sum_{\ell=1}^{2^k}  \Vert \E_{-\ell } (f_{[2^{\alpha j}]} ) \Vert_p  +     \sum_{\ell=0}^{2^k-1}  \Vert f_{[2^{\alpha j}]}- {\mathbb E}_{\ell} (f_{[2^{\alpha j}]} )  \Vert_p \\
 + 2^{k/2} \sum_{\ell \geq  2^k +1} \Vert \E_{-\ell}(f_{[2^{\alpha j}]}  ) \Vert_p   + 2^{k/2}
 \sum_{\ell \geq  2^k } \Vert f_{[2^{\alpha j}]}  - \E_{\ell}(f_{[2^{\alpha j}]} ) \Vert_p \, .
\end{multline*}
Therefore, taking into account the upper bound \eqref{auto3} for the two first terms in the right hand side, and the upper bound \eqref{auto1} to handle the two last terms
(since $[2^{\alpha j}] \leq b(2^k)$), we derive that
\begin{equation} \label{3declma2}
\max_{1 \leq m \leq 2^k}\Vert R_{j,m} \Vert_p \ll    2^{k ( 1 - \theta (p-1)/p)}
\end{equation}
(recall that $\theta(p-1)/p<1$).
On the other hand, starting from (\ref{2declma2}) with $N=2 [2^{kp/2}]$ and using Lemma 5.1 in \cite{DMP1}, we get that
\begin{multline*}
\max_{1 \leq m \leq 2^k}\Vert R_{j,m} \Vert_p \ll  \sum_{\ell=1}^{2 [2^{kp/2}]}  \Vert \E_{-\ell } (f_{[2^{\alpha j}]} ) \Vert_p  +     \sum_{\ell=0}^{2 [2^{kp/2}]}  \Vert f_{[2^{\alpha j}]}- {\mathbb E}_{\ell} (f_{[2^{\alpha j}]} )  \Vert_p \\
 + 2^{k/2} \sum_{\ell \geq  [2^{kp/2}]} \frac{\Vert \E_{-\ell}(f_{[2^{\alpha j}]}  ) \Vert_p  }{\ell^{1/p}} + 2^{k/2}
 \sum_{\ell \geq  [2^{kp/2}] } \frac{ \Vert f_{[2^{\alpha j}]}  - \E_{\ell}(f_{[2^{\alpha j}]} ) \Vert_p }{\ell^{1/p}}\, .
\end{multline*}
Therefore, it follows from \eqref{auto3} that
\begin{equation} \label{3declma2*}
\displaystyle \max_{1 \leq m \leq 2^k}\Vert R_{j,m} \Vert_p  \ll    2^{(kp ( 1 - \theta (p-1)/p))/2} \, .
\end{equation}
Let
\beq \label{defCj0}
C = [\alpha^{-1 }( \log ( \gamma^{-1}))/(\log 2)] \, \text{ and } j_0= (\log 2)^{-1}(\log j - \log C ) \, .
\eeq  Clearly, if $j_0 \leq  k $ then $[2^{\alpha j}] \leq b(2^k)$. Therefore using the  upper bound \eqref{3declma2*} when $k < j_0$ and the upper bound \eqref{3declma2} when $k \geq j_0$, we get that for any positive integer $j$
$$
\sum_{k=0}^{j} 2^{-k/p} \Vert R_{j,2^k} \Vert_p \ll \frac{j^{p/2}}{j^{1/p}  j^{\theta (p-1)/2}}\, ,
$$
since $\theta<(p^2-2)/(p(p-1))$.
Now, since $\theta >1$, it follows that
$$
\sum_{j\ge 1} j^{-p}\Big ( \sum_{k=0}^{j} 2^{-k/p} \Vert R_{j,2^k} \Vert_p \Big ) ^p  <\infty \, .
$$
From \eqref{1declma2}, this implies that \eqref{lma1p1} holds with $D_j$ defined by \eqref{0declma2} and $q=p$. This ends the proof of lemma \ref{Lem2}. $\square$

\bigskip

\noindent{\bf Proof of Lemma \ref{Lem3}.}  Let $M_n = \sum_{\ell=1}^n d_{\ell} $ where $d_{\ell}= d_0 \circ T^{\ell}$ with $d_0 = \sum_{i \in {\mathbb Z}} P_0 (f \circ T^i)$. Notice that the upper bound (\ref{auto4}) and the fact that $\theta >1$ imply in particular that $$
\sum_{n >1 } n^{-1/p}\Vert \E_{-n} ( f  ) \Vert_p < \infty \text{ and } \sum_{n >1 } n^{-1/p}\Vert f - \E_n ( f )  \Vert_p < \infty \, ,
$$
and then that $\sum_{k \in {\mathbb Z}} \Vert P_0 (X_k) \Vert_p < \infty$ (use for instance Lemma 5.1 in \cite{DMP1} to see this). Therefore $\sum_{k \in {\mathbb Z}} \Vert P_0 (X_k) \Vert_2 < \infty$. Using (\ref{app1}) of Lemma \ref{approxrnp}, we get that
\begin{equation} \label{app2**}
\Vert S_n (f) - M_n \Vert_2 = o (\sqrt{n}) \, .
\end{equation}
From
(\ref{defsigma2})
of Theorem \ref{asip}, we know that $ n^{-1} \Vert S_n(f) \Vert_2^2$
converges to $\sigma^2$.  It follows from \eqref{app2**}  that
$\sigma^2 =n^{-1} \E (M_n^2)= \E (d_0^2)$.

We shall prove now that
\begin{equation} \label{condshao3}
\Vert M_n^*- M_n \Vert_2 =O\big ( n^{1/p} \big ) \, .
\end{equation}
Let $N$ be the positive integer such that $ 2^{N-1} < n  \leq 2^N$. Since  $M_n^*- M_n$,
is a martingale, we have that
\begin{equation} \label{estnorm2}
\Vert M_n^*- M_n \Vert^2_2 = \sum_{\ell = 1}^n \E ( (d_{\ell}^* - d_{\ell})^2) \leq \E ( (d_{1}^* - d_{1})^2) +  \sum_{j=0}^{N-1} \sum_{\ell = 2^{j}+1}^{2^{j+1}} \E ( (d_{\ell}^* - d_{\ell})^2) \, .
\end{equation}
By stationarity, for any $\ell \in [2^{j}+1, 2^{j+1}] \cap {\mathbb N}$ we get that
\begin{multline*}
\Vert d_{\ell}^* - d_{\ell}\Vert_2 = \Vert \sum_{i \in {\mathbb Z}} P_0 ( ( f -f_{[2^{\alpha j}]}) \circ T^i) \Vert_2 \\
\leq 2^{j+3} \Vert  f -f_{[2^{\alpha j}]}  \Vert_2 +  \sum_{i \geq 2^{j+1}}\Vert  P_{-i} (  f -f_{[2^{\alpha j}]} ) \Vert_2 + \sum_{i \geq 2^{j+1}}\Vert  P_{i} (  f -f_{[2^{\alpha j}]}) \Vert_2 \, .
\end{multline*}
According to \eqref{majdiff2}
 \begin{equation} \label{est1lma3}
\Vert f -f_{[2^{\alpha j}]}  \Vert_2\le R 2^{-\zeta \alpha j/ 2} \, .
\end{equation}
On the other hand, by Lemma 5.1 in \cite{DMP1},
$$
 \sum_{i \geq 2^{j+1}} \Vert  P_{-i} (  f -f_{[2^{\alpha j}]} ) \Vert_2 \ll \sum_{k \geq 2^{j}}  k^{-1/2}\Vert  \E_{-k} (  f -f_{[2^{\alpha j}]} ) \Vert_2 $$
 and
$$ \sum_{i \geq 2^{j+1}}\Vert  P_{i} (  f -f_{[2^{\alpha j}]}) \Vert_2 \ll \sum_{k \geq 2^{j}} k^{-1/2}\Vert   f -f_{[2^{\alpha j}]}  +\E_{k} (  f -f_{[2^{\alpha j}]} ) \Vert_2 \, .
 $$
 Using the estimate \eqref{auto6} and \eqref{auto6bis}, it follows that 
 \begin{equation} \label{est2lma3}
 \sum_{i \geq 2^{j+1}} \big ( \Vert  P_{-i} (  f -f_{[2^{\alpha j}]} ) \Vert_2 +
\Vert  P_{-i} (  f -f_{[2^{\alpha j}]} ) \Vert_2 \big ) = O(\beta^{2^j})
 \, .
  \end{equation}
Combining  the upper bounds \eqref{est1lma3} and \eqref{est2lma3}  with the fact that $\alpha \zeta \geq 3 -2/p   $, it follows that
  $$
  \Vert d_{\ell}^* - d_{\ell}\Vert_2 \ll  2^{-j (p-2)/(2p) } \, .
  $$
Using this estimate in \eqref{estnorm2}, we obtain that
$
\Vert M_n^*- M_n \Vert^2_2 \ll n^{2/p}
$,
 proving \eqref{condshao3}.

Now, let us recall Theorem 2.1 in \cite{shao93} (used with $a_n = n^{2/p} ( \log n)$):
if there exists a finite constant $K$ such that
\begin{equation} \label{condshao1}
\sup_{k \geq 1} \Vert d^*_k \Vert_p \leq K \, ,
\end{equation}
and if
\begin{equation} \label{condshao2}
\sum_{i=1}^n\left({\mathbb E}((d^*_i)^2|\mathcal F_{i-1})-{\mathbb E}((d^*_i)^2)\right)
     =o\big ( n^{2/p} (\log n)\big )\ \ a.s. \, ,
\end{equation}
then, since $\E ((M^*_n)^2)\sim n\sigma^2$, enlarging ${\mathbb T}^d$ if
necessary, there exists a sequence $(Z^*_\ell)_{\ell
\geq 1}$ of independent  Gaussian random variables with zero mean and
variance $\E( Z_\ell^*)^2= \E (d_{\ell}^*)^2= (\sigma_{\ell}^*)^2$ such that
\begin{equation} \label{asipmart1}
\sup_{1\leq k \leq n} \Big| M^*_k - \sum_{\ell=1}^k Z^*_{\ell}\Big|  = o \big ( n^{1/p} ( \log n)  \big )
\text{ almost surely, as $n\rightarrow \infty$}.
\end{equation}
Let $(\delta_k)_{k \geq
1}$ be a sequence of iid \ Gaussian random variables with mean zero
and variance $\sigma^2$, independent of the sequence $(Z^*_{\ell})_{\ell\geq
1}$. We now construct a sequence $(Z_{\ell})_{\ell \geq 1}$ as follows. If $\sigma_{\ell}^*=0$, then $Z_\ell=\delta_\ell$,
else $Z_\ell=
(\sigma/\sigma_{\ell}^*) Z_{\ell}^*$. By construction, the $Z_\ell$'s are iid \ Gaussian random variables with
mean zero and variance $\sigma^2$. Let $G_\ell=Z_\ell-Z^*_{\ell}$ and note that $(G_\ell)_{\ell \geq 1}$ is a
sequence of independent Gaussian random variables with mean zero and
variances $\Var(G_\ell)=(\sigma-\sigma_{\ell}^*)^2$. Notice now that
\[
  v^2_n= \Var\Big( \sum_{i=1}^n G_i \Big) = \sum_{i=1}^n \big ( \Vert d_i \Vert_2 -\Vert d^*_i \Vert_2 \big )^2\leq  \Vert M_n -M_n^* \Vert^2_2 \, .
\]
From the basic inequality
\[
  {\mathbb P}\Big
  ( \max_{1 \leq k \leq n} \Big| \sum_{i=1}^k G_i \Big|>x\Big) \leq
  2 \exp \Big(-\frac{ x^2}{2  v^2_n}\Big)\, ,
\]
and the fact that by \eqref{condshao3}, $v_n^2\ll n^{2/p}$, it follows that for any $\varepsilon >0$,
  \begin{equation*}
   \sum_{n >1} n^{-1}{\mathbb P}\Big
  ( \max_{1 \leq k \leq n} \Big| \sum_{i=1}^k G_i \Big|> \varepsilon n^{1/p} ( \log n)\Big) < \infty \, ,
  \end{equation*}
showing that
$\max_{1 \leq k \leq n} | \sum_{i=1}^k G_i |= o  ( n^{1/p} ( \log n) )$ almost surely. Therefore starting from \eqref{asipmart1}, we conclude
that if \eqref{condshao1} and \eqref{condshao2} hold then Lemma \ref{Lem3} does; namely, enlarging ${\mathbb T}^d$ if
necessary, there exists a sequence $(Z_\ell)_{\ell
\geq 1}$ of iid  Gaussian random variables with zero mean and
variance $\sigma^2$ such that
\begin{equation*} \label{asipmart2}
\sup_{1\leq k \leq n} \Big| M^*_k - \sum_{\ell=1}^k Z_{\ell}\Big|
= o (n^{1/p} ( \log n)) \quad
\text{almost surely, as $n\rightarrow \infty$}.
\end{equation*}
It remains to show that \eqref{condshao1} and \eqref{condshao2} are satisfied. We start with \eqref{condshao1}. Notice that $\Vert d_1^* \Vert_p \leq \sum_{k\in\mathbb Z} \Vert P_{-k}\left(f_{1} \right)\Vert_p $ and that, for every $j\ge 0$ and every $\ell\in\{2^{j}+1,...,2^{j+1}\}$,
$$ \Vert d^*_\ell \Vert_p \leq \sum_{k\in\mathbb Z} \Vert P_{-k} (f_{[2^{\alpha j}]})  \Vert_p \, .$$
By Lemma 5.1 in \cite{DMP1},
$$
 \sum_{k\in\mathbb Z} \Vert P_{-k} (g )  \Vert_p \ll  \sum_{k \geq 1} k^{-1/p}\Vert \E_{-k} (g)
 \Vert_p +\sum_{k \geq 1 } k^{-1/p}\Vert g - \E_k ( g )  \Vert_p \, ,
$$
with constants non depending on $g$. Hence,
using the estimate \eqref{auto3} and the fact that $\theta >1$, we get that,
for every $j\ge 0$ and every $\ell\in\{2^{j}+1,...,2^{j+1}\}$, there exists a constant $K$ non depending on $j$ such that $\Vert d_{\ell}^* \Vert_p \leq K$. This ends the proof of \eqref{condshao1}.

\medskip

To prove (\ref{condshao2}), we proceed as follows. Following the beginning of the proof
of Lemma \ref{Lem1}, we infer that (\ref{condshao2}) will be proven if we can show that
\begin{eqnarray} \label{0declma3}
D_j := \sup_{ 1 \leq \ell \leq 2^{j}} \Big |  \sum_{i = 2^j +1}^{2^j + \ell} ({\mathbb E}((d^*_i)^2|\mathcal F_{i-1})-{\mathbb E}((d^*_i)^2)) \Big |
 = o\big (j \, 2^{2j/p}\big ) \ \ a.s. \, .
\end{eqnarray}
This will hold true as soon as
\begin{equation} \label{lma3p1}
\sum_{j\ge 1}\frac{ \Vert D_j \Vert_{p/2}^{p/2 }}{2^j\, j^{ p/2}}<\infty \, .
\end{equation}
For any $j$ fixed and any $i \in {\mathbb Z}$, let $d_{j,i} = \sum_{k\in\mathbb Z}  P_{i}\left(f_{[2^{\alpha j}]} \circ T^k \right)$. By stationarity
$$
\|D_j\|_{p/2} :=
\Big \|\sup_{ 1 \leq \ell \leq 2^{j}} \Big |  \sum_{i = 1}^{ \ell} ({\mathbb E}(d_{j,i}^2|\mathcal F_{i-1})-{\mathbb E}(d_{j,i}^2)) \Big | \Big \|_{p/2}\, .
$$
Observe now that, for any $j$ fixed, $(d_{j,i})_{i \in {\mathbb Z}}$ is a stationary sequence of martingale differences in ${\mathbb L}^p$. Let
$$
M_{j,n}:= \sum_{i=1}^n d_{j,i} \, .
$$
Applying Theorem 3 in \cite{WZ} (since $1 < p/2 \leq 2$) and using the martingale property of the sequence $(M_{j,n})_{n \geq 1}$, we get that
$$
\E \Big (  \sup_{ 1 \leq \ell \leq 2^{j}} \big |  \sum_{i = 1}^{ \ell} ({\mathbb E}(d_{j,i}^2|\mathcal F_{i-1})-{\mathbb E}(d_{j,i}^2)) \big |  \Big |^{p/2}  \Big ) \ll  2^j \| d^2_{j,1} \|^{p/2}_{p/2}  + 2^j \left ( \sum_{k=0}^{j-1} \frac{\|\E(M_{j,2^k}^2|{\mathcal F}_{0}) - \E(M_{j,2^k}^2) \|_{p/2}}{2^{2k/p}}  \right )^{{p/2}} \, .
$$
Using the fact that $\| d^2_{j,1} \|_{p/2} = \| d_{j,1} \|^2_{p} \leq K$ where $K$ does not depend on $j$, the convergence \eqref{lma3p1} will be then proven if we can show that
$$
\sum_{j\ge 1} \frac{1}{j^{p/2}} \Big ( \sum_{k=0}^{j-1} 2^{-2k/p} \|\E_0(M_{j,2^k}^2) - \E(M_{j,2^k}^2) \|_{p/2} \Big )^{{p/2}}<\infty \, .
$$
According to the arguments developed in the proof of Theorems 3.1 and 3.2 in \cite{DMP1}
(see (3.19) and (3.20) of \cite{DMP1}),
since $ (M_{j,k})_{k \geq 1}$ is a sequence of martingales, we infer that this last convergence will be satisfied as soon as there exists a positive integer $c$ such that
\beq \label{condcarremart}
\sum_{j\ge 1}\frac{ 1}{ j^{p/2}} \Big (  \sum_{k=0}^j  2^{-2k/p}\big \Vert  \E_{-c2^k } (M_{j,2^k}^2  ) -\E  (M_{j,2^k}^2  )
 \big \Vert_{p/2} \Big )^{p/2}< \infty \, .\eeq We shall prove in what follows that this convergence holds as soon as $c$ is chosen in such a way that \eqref{auto5} holds true.

For any positive integer $n$, let
$$
S_{j,n} = \sum_{\ell=1}^n f_{[2^{\alpha j}]} \circ T^{\ell} \, \text{ and } \, R_{j,n} = S_{j,n} -M_{j,n} \, .
$$
We first write that
\begin{multline} \label{condcarremartP0}
\Vert  \E_{-c2^k } (M_{j,2^k}^2  ) -\bkE  (M_{j,2^k}^2  )
 \Vert_{p/2} \leq \Vert  \E_{-c2^k } (S_{j,2^k}^2  ) -\bkE  (S_{j,2^k}^2  )
  \Vert_{p/2} \\
 +4 \Vert \E_{-c2^k } (S_{j,2^k}R_{j,2^k}  ) \Vert_{p/2} +
 2 \Vert R_{j,2^k} \Vert^{2}_{p} \, .
\end{multline}
Let $j_0$ be defined as in \eqref{defCj0}. Using the  upper bound \eqref{3declma2*} when $k < j_0$ and the upper bound \eqref{3declma2} when $k \geq j_0$,
we get that for any positive integer $j$,
$$
\sum_{k=0}^{j} 2^{-2k/p} \max_{1 \leq m \leq 2^k}\Vert R_{j,m} \Vert^2_p \ll
\frac{j^{p}}{j^{2/p}  j^{\theta (p-1)}} \, ,
$$
since we can assume without loss of generality that $\theta<(p^2-1)/(p(p-1))$. 
Now, since $\theta >1$, it follows that
\beq \label{condcarremartP1}
\sum_{j\ge 1}\frac{ 1}{ j^{ p/2}} \Big (  \sum_{k=0}^j  2^{-2k/p}\max_{1 \leq m \leq 2^k}\Vert R_{j,m} \Vert^2_p  \Big )^{p/2}< \infty \, .\eeq
On an other hand, $c$ being chosen such that \eqref{auto5} holds true, the upper bound in \eqref{auto5} together with the fact that $\theta >1$, implies that
$$
\sum_{k=0}^{j} 2^{-2k/p} \max_{1 \leq m \leq 2^k} \Vert  \E_{-c2^k } (S_{j,m}^2  ) -\bkE  (S_{j,m}^2  )
 \big \Vert_{p/2} \ll \sum_{k=0}^{j} 2^{-2k/p}  2^{2k (1 - \theta (p-1)/p)} = O(1) \, .
$$
Therefore,
\beq \label{condcarremartP2}
\sum_{j\ge 1}\frac{ 1}{ j^{ p/2}} \Big (  \sum_{k=0}^j  2^{-2k/p} \max_{1 \leq m \leq 2^k} \Vert  \E_{-c2^k } (S_{j,m}^2  ) -\bkE  (S_{j,m}^2  )
 \big \Vert_{p/2}   \Big )^{p/2}< \infty \, .\eeq
 Starting from \eqref{condcarremartP0}, and taking into account \eqref{condcarremartP1} and \eqref{condcarremartP2}, we then infer that \eqref{condcarremart} will hold true if we can show that
 \beq \label{condcarremartP3}
\sum_{j\ge 1}\frac{ 1}{ j^{ p/2}} \Big (  \sum_{k=0}^j  2^{-2k/p}\Vert  \E_{-c2^k } (S_{j,2^k}R_{j,2^k}  )
 \big \Vert_{p/2}  \Big )^{p/2}< \infty \, .\eeq
With this aim, we use Inequality (3.24) in \cite{DMP1} (taking $n=2^k$,
$u_n=[2^{k/2}]$, $r=c2^k$). Therefore,
\begin{multline} \label{doubleprod}
\Vert  \E_{-c2^k } (S_{j,2^k}R_{j,2^k}  )
 \big \Vert_{p/2} \ll 2^{k/4} \big ( \Vert \E_{0}  (S_{j,2^k})\Vert_{2} +
 \Vert S_{j,2^k} - \E_{2^k} ( S_{j,2^k}) \Vert_{2} \big ) \\ + \max_{m=\{2^k,2^{k} - [2^{k/2}]\}} \Vert R_{j,m}\Vert_p^2 +
  2^{k/2} \big ( \Vert \E_{-[2^{k/2}]}  (S_{j,2^k} )\Vert_{2} + \Vert S_{j,2^k} - \E_{2^k + [2^{k/2}]}  (S_{j,2^k}  )\Vert_{2} \big )  \\ + \max_{m=\{2^k, [2^{k/2}]\}} \Vert \E_{-c2^k}  (S^2_{j,m} ) - \E(S_{j,m}^2)\Vert_{p/2}
 + 2^k  \sum_{|\ell| \geq 2^k} \Vert  P_0( f_{[2^{j \alpha}]}\circ T^{\ell}) \Vert_2  \, .
\end{multline}
By stationarity,
$$
\Vert \E_{0}  (S_{j,2^k})\Vert_{2} +
 \Vert S_{j,2^k} - \E_{2^k} ( S_{j,2^k}) \Vert_{2} \leq \sum_{\ell=1}^{2^k} \Vert \E_{-{\ell}}  (f_{[2^{j \alpha}]})\Vert_{2} + \sum_{\ell=0}^{2^k-1}
 \Vert f_{[2^{j \alpha}]} - \E_{\ell} ( f_{[2^{j \alpha}]}) \Vert_{2} \, .
$$
Hence by using \eqref{auto6},
\beq \label{p1doubleprod}
2^{k/4} \big ( \Vert \E_{0}  (S_{j,2^k})\Vert_{2} +
 \Vert S_{j,2^k} - \E_{2^k} ( S_{j,2^k}) \Vert_{2} \big ) \ll 2^{k/4} \, .
\eeq
Using again
\eqref{auto6}
and the
stationarity, we get that
\begin{align*}
 \Vert \E_{-[2^{k/2}]} &  (S_{j,2^k} )\Vert_{2} + \Vert S_{j,2^k} - \E_{2^k + [2^{k/2}]}  (S_{j,2^k}  )\Vert_{2} \\
 &\leq \sum_{\ell=1}^{2^k} \Vert \E_{-([2^{k/2}] + {\ell})}  (f_{[2^{j \alpha}]})\Vert_{2} + \sum_{\ell=0}^{2^k-1}
 \Vert f_{[2^{j \alpha}]} - \E_{[2^{k/2}] + \ell} ( f_{[2^{j \alpha}]}) \Vert_{2}
  \ll \beta^{2^{k/2}} \, .
\end{align*}
Therefore
\beq \label{p2doubleprod}
2^{k/2} \big ( \Vert \E_{-[2^{k/2}]}  (S_{j,2^k} )\Vert_{2} + \Vert S_{j,2^k} - \E_{2^k + [2^{k/2}]}  (S_{j,2^k}  )\Vert_{2} \big ) =O(1) \, .
\eeq
On an other hand, using again
\eqref{auto6} and the
stationarity,
\begin{align*}
 \sum_{|\ell| \geq 2^k} \Vert  P_0( f_{[2^{j \alpha}]}\circ T^{\ell}) \Vert_2 & \leq \sum_{\ell \geq 2^k}  \big ( \Vert  \E_{- \ell}( f_{[2^{j \alpha}]}) \Vert_2 +
\Vert  f_{[2^{j \alpha}]} - \E_{ \ell -1}( f_{[2^{j \alpha}]}) \Vert_2 \big )
 \ll  \beta^{2^{k}} \, ,
\end{align*}
which implies that
\beq \label{p3doubleprod}
2^k  \sum_{|\ell| \geq 2^k} \Vert  P_0( f_{[2^{j \alpha}]}\circ T^{\ell}) \Vert_2 = O (1) \, .
\eeq
Starting from \eqref{doubleprod} and taking into account the convergence \eqref{condcarremartP1} and \eqref{condcarremartP2}, and the upper bounds \eqref{p1doubleprod}, \eqref{p2doubleprod} and \eqref{p3doubleprod}, we then derive that \eqref{condcarremartP3} holds. This ends the proof of \eqref{condshao2} and therefore of Lemma \ref{Lem3}. $\square$

\section{Appendix}
\setcounter{equation}{0}
Let $(\Omega,{\cal A}, {\mathbb P} )$ be a
probability space, and $T:\Omega \mapsto \Omega$ be
 a bijective bimeasurable transformation preserving the probability ${\mathbb P}$. Let us denote by $| \cdot |_{m}$ the euclidean norm on ${\mathbb R}^m$ and by $< \cdot, \cdot >_m$ the associated scalar product. For a $\sigma$-algebra ${\mathcal F}_0$ satisfying ${\mathcal F}_0
\subseteq T^{-1 }({\mathcal F}_0)$, we define the nondecreasing
filtration $({\mathcal F}_i)_{i \in {\mathbb Z}}$ by ${\mathcal F}_i =T^{-i
}({\mathcal F}_0)$. Let ${\mathcal {F}}_{-\infty} = \bigcap_{k \in {\mathbb
Z}} {\mathcal {F}}_{k}$ and ${\mathcal {F}}_{\infty} = \bigvee_{k \in
{\mathbb Z}} {\mathcal {F}}_{k}$.
For a random variable $X$ with values in ${\mathbb R}^m$, we denote   by
$\Vert X \Vert_{p,m}= (  \E ( | X|^p _{m}))^{1/p}$
its norm in ${\mathbb L}^{p} ({\mathbb R}^m)$.

In what follows $X_0$ is a random variable with values in ${\mathbb R}^m$, and we define the stationary sequence  $(X_i)_{i \in \mathbb Z}$ by
$X_i = X_0 \circ T^i$. We shall use the notations $\E_k (X) = \E (X | {\mathcal F}_k)$, $\E_\infty (X) = \E (X | {\mathcal F}_\infty)$,
 $\E_{-\infty} (X) = \E (X | {\mathcal F}_{-\infty})$, and $P_k(X)= \E_{k}(X)-\E_{k-1}(X)$.

The aim of this section is to collect some results about invariance principles for stationary sequences that are non necessarily adapted to the underlying filtration.
We start with a martingale approximation result. The estimate \eqref{app1}
of Proposition \ref{Gor} below is a generalization of Item 2 of Theorem 1 in \cite{W} to the multidimensional case and to the case where the variables are non necessarily adapted to the filtration under consideration. The convergence \eqref{app2} is new. Notice that the proof of the next lemma is based on algebraic computations and on Burkholder's inequality.  Burkholder's inequality being also valid in Hilbert spaces (see \cite{BU}), the approximation lemma below is then also valid for variables taking values in a  separable Hilbert space, ${\mathcal H}$, by replacing the norm $| \cdot |_{m}$ by the norm on ${\mathcal H}$, let say $| \cdot |_{{\mathcal H}}$.

\begin{Proposition} \label{approxrnp} Let  $p\in [1, \infty[$ and $p'=\min ( 2,p)$.
Assume that $\E_{- \infty} (X_0) = 0$ almost surely, that $\E_{\infty} (X_0) = X_0$ almost surely, and that
\begin{equation}\label{Gor}
  \sum_{i \in {\mathbb Z}} \Vert P_0(X_i) \Vert_{p,m} < \infty \, .
\end{equation}
Let $d_0=\sum_{i \in {\mathbb Z}} P_0(X_i)$, $M_n:= \sum_{i=1}^n d_0 \circ T^i $ and $R_n:=S_n-M_n$. For any positive integer $n$,
\beq \label{app1}
\| R_n \|^{p'}_{p,m}  \ll \sum_{k= 1}^{n} \Big ( \sum_{|\ell| \ge k}\| P_{\ell}(X_0)\|_{p,m} \Big )^{p'} \, .
\eeq
In addition,
\beq \label{app2}
\Big \| \max_{1 \leq k \leq n} | R_k |_m
\Big \|_{p}  = o(n^{1/p'}) \ \text{  as $n\rightarrow \infty$}.
\eeq
\end{Proposition}
\begin{Remark}
The constant
appearing in (\ref{app1}) depends only on $p$ and not on $(\Omega,\mathcal A,{\mathbb P},T,X_0,
\mathcal F_0)$.
\end{Remark}

\noindent{\bf Proof of Proposition \ref{approxrnp}.}
It will be useful to note that $P_i(X_j)=P_{i-j}(X_{0})\circ T^j=P_{0}(X_{j-i})\circ T^i$ almost surely.
The following decomposition is valid:
\begin{eqnarray} \label{decRn}
R_n   &= & \sum_{k=1}^n \Big ( X_k - \sum_{j =1}^n P_j(X_k) \Big )
  - \sum_{k=1}^n \sum_{j \geq n+1} P_k(X_j) - \sum_{k=1}^n \sum_{j=0}^{\infty} P_k(X_{-j}) \nonumber \\
  &= & {\mathbb E}_0 (S_n) - \sum_{k=1}^n \sum_{j \geq n+1} P_k(X_j) + S_n - {\mathbb E}_n (S_n )  -
\sum_{k=1}^n \sum_{j=0}^{\infty} P_k(X_{-j}) \, .
\end{eqnarray}
Applying Burkholder's inequality for multivariate martingales, and using
the  stationarity, we obtain that
there exists a positive constant $c_p$ such that, for any positive integer $n$,
\begin{eqnarray} \label{Burkh1}
\Big \Vert  \sum_{k=1}^n \sum_{j \geq n+1} P_k(X_j) \Big \Vert^{p'}_{p,m}   \leq   c_p
  \sum_{k=1}^n \Big \Vert \sum_{j \geq n +1} P_k(X_j) \Big \Vert_{p,m}^{p'}  \leq  c_p
   \sum_{k=1}^n  \Big (   \sum_{j \geq k} \Vert P_0(X_j)  \Vert_{p,m} \Big ) ^{p'} \, ,
\end{eqnarray}
and
\begin{eqnarray} \label{Burkh2}
\Big \Vert  \sum_{k=1}^n \sum_{j\geq 0} P_k(X_{-j})  \Big \Vert^{p'}_{p,m}    \leq
c_p  \sum_{k=1}^n \Big \Vert \sum_{j\geq 0} P_k(X_{-j})  \Big \Vert_{p,m} ^{p'}  = c_p
\sum_{k=1}^n  \Big (   \sum_{j \geq k}  \Vert P_0(X_{-j}) \Vert_{p,m} \Big )^{p'} \, .
\end{eqnarray}

On an other hand, since $\E_{-\infty}(X_0)=0$ almost surely,
we have $\E_0(S_n)=\sum_{k\ge 0}
P_{-k}(S_n)$ almost surely. Hence by Burkholder's inequality for multivariate martingales together with stationarity, there exists a positive constant $c_p$ depending only on $p$ such that
\begin{multline*}
\|\E_0(S_n)\|_{p,m}^{p'}\le c_p \sum_{k\ge 0}\|P _{-k}(S_n)\|_{p,m}^{p'} \leq c_p \sum_{k\ge 0} \Big ( \sum_{\ell=1}^n\|P _{-k}(X_{\ell})\|_{p,m} \Big )^{p'} \\ \le
c_p\sum_{k= 0}^{n-1} \Big ( \sum_{\ell \ge k+1}\| P_{-\ell}(X_0)\|_{p,m} \Big )^{p'}
+c_p \sum_{k\ge n} \Big (\sum_{\ell=k+1}^{k+n}\| P_{-\ell}(X_0)\|_{p,m} \Big )^{p'}\\
\le c_p\sum_{k= 0}^{n-1} \Big ( \sum_{\ell \ge k+1}\| P_{-\ell}(X_0)\|_{p,m} \Big )^{p'}
+c_p \Big (\sum_{i \ge n+1}\|P_{-i}(X_0)\|_{p,m} \Big )^{p'-1}
\sum_{k\ge n} \sum_{\ell=k+1}^{k+n}\| P_{-\ell}(X_0)\|_{p,m} \\
\le c_p\sum_{k= 0}^{n-1} \Big ( \sum_{\ell \ge k+1}\| P_{-\ell}(X_0)\|_{p,m} \Big )^{p'} + c_p
n \Big (\sum_{\ell \ge n+1}\|P_{-\ell}(X_0)\|_{p,m} \Big )^{p'}  \, .
\end{multline*}
Therefore
\beq \label{majeosn}
\|\E_0(S_n)\|_{p,m}^{p'} \leq 2 c_p\sum_{k= 1}^{n} \Big ( \sum_{\ell \ge k}\| P_{-\ell}(X_0)\|_{p,m} \Big )^{p'}  \, .
\eeq
We handle now the quantity $\Vert S_n - {\mathbb E}_n (S_n ) \Vert_{p,m}^{p'} $. Since
 ${\mathbb E}_\infty(X_0)=X_0$ almost surely, we first write that $S_n - {\mathbb E}_n (S_n ) = \sum_{k \geq n+1} P_k (S_n)$. Hence,
applying Burkholder's inequality for multivariate martingales and using the  stationarity, we infer that there exists a positive constant $c_p$ depending only on $p$ such that
\begin{multline*}
\|S_n - {\mathbb E}_n (S_n )\|_{p,m}^{p'}\le c_p \sum_{k\ge n+1}\|P _{k}(S_n)\|_{p,m}^{p'} \leq c_p \sum_{k\ge n+1} \Big ( \sum_{\ell=1}^n\|P _{k}(X_{\ell})\|_{p,m} \Big )^{p'} \\ \le
c_p\sum_{k= n+1}^{2n} \Big ( \sum_{\ell=k-n}^{k-1}\| P_{\ell}(X_0)\|_{p,m} \Big )^{p'}
+c_p \sum_{k\ge 2n +1} \Big (\sum_{\ell=k-n}^{k-1}\| P_{\ell}(X_0)\|_{p,m} \Big )^{p'}\\
\le c_p\sum_{k= 1}^{n} \Big ( \sum_{\ell \ge k}\| P_{\ell}(X_0)\|_{p,m} \Big )^{p'}
+c_p \Big (\sum_{i \ge n+1}\|P_{i}(X_0)\|_{p,m} \Big )^{p'-1}
\sum_{k\ge 2n +1} \sum_{\ell=k-n}^{k-1}\| P_{\ell}(X_0)\|_{p,m}  \\
\le c_p\sum_{k= 1}^{n} \Big ( \sum_{\ell \ge k}\| P_{\ell}(X_0)\|_{p,m} \Big )^{p'} + c_p
n \Big (\sum_{\ell \ge n+1}\|P_{\ell}(X_0)\|_{p,m} \Big )^{p'}  \, .
\end{multline*}
Therefore
\beq \label{majeosnnonadap}
\|S_n - \E_n (S_n) \|_{p,m}^{p'} \leq 2 c_p\sum_{k= 1}^{n} \Big ( \sum_{\ell \ge k}\| P_{\ell}(X_0)\|_{p,m} \Big )^{p'}  \, .
\eeq
Starting from \eqref{decRn} and taking into account the upper bounds \eqref{Burkh1}, \eqref{Burkh2}, \eqref{majeosn} and \eqref{majeosnnonadap}, the inequality \eqref{app1} follows.

\medskip

We turn now to the proof of \eqref{app2}. 
Let $r$ be some fixed positive integer. 
Since $M_k=\sum_{i=1}^kd_0\circ T^i$ and  $d_0\circ T^i=\sum_{j\in\mathbb Z}P_i(X_j)$,
 the following decomposition holds:
\begin{gather} \label{decRnbis}
R_k   =  \sum_{i=1}^k X_i -  \sum_{i =1}^k \sum_{\ell = -r-i}^{r+i} P_i(X_{\ell})
  - \sum_{i =1}^k \sum_{\ell \geq r+i +1 } P_i(X_{\ell})-\sum_{i =1}^k \sum_{\ell \geq r+i +1 }
P_i(X_{- \ell}) \, .
\end{gather}
Applying Burkholder's inequality for multivariate martingales and using the stationarity, we infer that there exists a positive constant $c_p$ depending only on $p$ such that, for any positive integer $n$,
\begin{eqnarray} \label{Burkh1bis}
\Big \Vert  \max_{r \leq k \leq n}  \Big | \sum_{i =1}^k \sum_{\ell \geq r+i +1 } P_i(X_{\ell}) \Big |_m   \Big \Vert^{p'}_{p}   \leq   c_p
\sum_{i =1}^n \Big \Vert \sum_{\ell \geq r+i +1 } P_i(X_{\ell}) \Big \Vert_{p,m}^{p'}  \leq  c_p
n \Big (   \sum_{j \geq r+1} \Vert P_0(X_j)  \Vert_{p,m} \Big ) ^{p'} ,
\end{eqnarray}
since $P_i(X_\ell)=P_0(X_{\ell-i})\circ T^i$.
Similarly
\begin{eqnarray} \label{Burkh2bis}
\Big \Vert  \max_{r \leq k \leq n}  \Big | \sum_{i =1}^k \sum_{\ell \geq r+i +1 } P_i(X_{-\ell}) \Big |_m   \Big \Vert^{p'}_{p}   \leq    c_p
n \Big (   \sum_{j \geq r+1} \Vert P_0(X_{-j})  \Vert_{p,m} \Big ) ^{p'} \, .
\end{eqnarray}
We write now that
\begin{gather} \label{decRnter}
\sum_{i=1}^k X_i -  \sum_{i =1}^k \sum_{\ell = -r-i}^{r+i} P_i(X_{\ell}) = \sum_{i=1}^k X_i -  \sum_{i =1}^k \sum_{\ell = 1}^{r+i} P_i(X_{\ell})-\sum_{i =1}^k \sum_{\ell = 0}^{r+i} P_i(X_{-\ell}) \, .
\end{gather}
The following decomposition holds
\begin{align}
  \sum_{i =1}^k \sum_{\ell = 1}^{r+i} P_i(X_{\ell})
   & =  \sum_{i =1}^{k-r} \sum_{\ell = 1}^{r+i}
P_i(X_{\ell}) + \sum_{i =k-r+1}^k \sum_{\ell = 1}^{r+i} P_i(X_{\ell})\nonumber  \\
   &=    \sum_{\ell = 1}^{k} \sum_{i =1}^{k-r} {\bf 1}_{i \geq \ell -r } P_i(X_{\ell})
   + \sum_{\ell = 1}^{k+r} \sum_{i =k-r+1}^k  {\bf 1}_{i \geq \ell -r } P_i(X_{\ell}) \, .
\label{decRn4}
\end{align}
Now,
\begin{align}
 \sum_{\ell = 1}^{k} \sum_{i =1}^{k-r} {\bf 1}_{i \geq \ell -r } P_i(X_{\ell})&  = \sum_{\ell = 1}^{r} \sum_{i =1}^{k-r}  P_i(X_{\ell}) +\sum_{\ell = r+1}^{k} \sum_{i =\ell -r }^{k-r} P_i(X_{\ell})\nonumber  \\
 &= \E_{k-r} (S_r) - \E_0 (S_r) + \E_{k-r} (S_k- S_r) - \sum_{\ell = r+1}^{k} \E_{\ell-r-1}
(X_{\ell})  \nonumber\\
&= \E_{k-r} (S_k) - \E_0 (S_r) - \sum_{\ell = r+1}^{k} \E_{\ell-r-1} (X_{\ell})\, ,\label{decRn5}
\end{align}
and
\begin{align}
 \sum_{\ell = 1}^{k+r} \sum_{i =k-r+1}^k  {\bf 1}_{i \geq \ell -r } P_i(X_{\ell}) & =
 \sum_{\ell = 1}^{k} \sum_{i =k-r+1}^k  P_i(X_{\ell}) +\sum_{\ell = k+1}^{k+r} \sum_{i =\ell -r }^{k} P_i(X_{\ell})
    \nonumber \\
& = \E_{k} (S_k) - \E_{k-r} (S_k) + \E_{k} (S_{k+r} - S_k) - \sum_{\ell = k+1}^{k+r}  \E_{\ell -r -1}(X_{\ell})
\nonumber\\
& = \E_{k} (S_{k+r}) - \E_{k-r} (S_k) - \sum_{\ell = k+1}^{k+r}  \E_{\ell -r -1}(X_{\ell})\, .
\label{decRn6}
\end{align}
Therefore starting from \eqref{decRnter}, and considering the decompositions \eqref{decRn4}, \eqref{decRn5} and \eqref{decRn6}, we get that
\beq \label{decRn7}
\sum_{i=1}^k X_i -  \sum_{i =1}^k \sum_{\ell = -r-i}^{r+i} P_i(X_{\ell})   = S_k
- \E_{k} (S_{k+r} ) + \E_0 (S_r) + \sum_{\ell = r+1}^{k+r} \E_{\ell-r-1} (X_{\ell})  -\sum_{i =1}^k \sum_{\ell = 0}^{r+i} P_i(X_{-\ell}) \, .
\eeq
The decomposition \eqref{decRnbis} together with the upper bounds \eqref{Burkh1bis}, \eqref{Burkh2bis} and \eqref{decRn7} imply that
\begin{multline} \label{decRn8}
\Big  \Vert  \max_{r \leq k \leq n}   | R_k |_m
\Big \Vert^{p'}_{p} \ll n \Big (   \sum_{|j| \geq r+1} \Vert P_0(X_j)  \Vert_{p,m} \Big ) ^{p'}
 + \big \Vert  \max_{r \leq k \leq n}  \big | S_k
- \E_{k} (S_{k+r} ) \big |_m \big \Vert^{p'}_{p} + \Vert \E_0 (S_r ) \Vert^{p'}_{p,m}\\ + \Big \Vert  \max_{r \leq k \leq n}  \Big | \sum_{\ell = r+1}^{k+r} \E_{\ell-r-1} (X_{\ell})  \Big |_m \Big \Vert^{p'}_{p}+\Big \Vert  \max_{r \leq k \leq n}  \Big | \sum_{i =1}^k \sum_{\ell = 0}^{r+i} P_i(X_{-\ell}) \Big | \Big \Vert^{p'}_{p}  \, .
\end{multline}
Applying Burkholder's inequality for multivariate martingales and using the  stationarity, there exists a positive constant $c_p$ depending only on $p$ such that, for any positive integer $n$,
\begin{eqnarray} \label{Burkh3bis}
\Big \Vert  \max_{r \leq k \leq n}  \Big | \sum_{i =1}^k \sum_{\ell = 0}^{r+i} P_i(X_{-\ell}) \Big | \Big \Vert^{p'}_{p}  \leq   c_p
\sum_{i =1}^n \Big \Vert \sum_{\ell = 0}^{r+i} P_i(X_{-\ell})  \Big \Vert_{p,m}^{p'}  \leq  c_p
\sum_{i =1}^n\Big (   \sum_{j \geq i} \Vert P_0(X_{-j})  \Vert_{p,m} \Big ) ^{p'} \, .
\end{eqnarray}
To handle the fourth term in the right-hand side of \eqref{decRn8} we proceed as follows. Since $\E_{- \infty} (X_\ell) = 0$ almost surely, we first write that
\begin{equation*}
\E_{\ell-r-1} (X_{\ell}) =\sum_{j=r+1 }^{\infty }P_{\ell-j}(X_{\ell}) \, .
\end{equation*}%
Then
\begin{equation*}
\max_{r \leq k \leq n}  \Big | \sum_{\ell = r+1}^{k+r} \E_{\ell-r-1} (X_{\ell})  \Big |_m \leq \sum_{j=r+1}^{\infty
} \max_{r \leq k \leq n}  \Big | \sum_{\ell = r+1}^{k+r}  P_{\ell-j}(X_{\ell}) \Big |_m \,.
\end{equation*}
Let now
$$
u_{i} = \Vert P_0(X_{i}) \Vert_{p,m} \, , \, C_r = \sum_{ i \geq r+1} u_{i} \, \text{ and } \, \alpha_{i } = C_r^{-1} u_{i} \, .
$$
By using the facts that for any $p \geq 1$, $x \mapsto x^p$ is convex and that $\alpha_{i}\geq
0$ with $\sum_{i \geq r+1} \alpha_{i}=1$ and writing $\sum_j a_j=\sum_j\alpha_j(a_j/\alpha_j)$,
we obtain that
$$
\Big \Vert  \max_{r \leq k \leq n}  \Big | \sum_{\ell = r+1}^{k+r} \E_{\ell-r-1} (X_{\ell})  \Big |_m \Big \Vert^{p}_{p}
\leq  \sum_{j=r+1}^{\infty
} \alpha_j^{1-p} \E \Big ( \max_{r \leq k \leq n}  \Big | \sum_{\ell = r+1}^{k+r}  P_{\ell-j}(X_{\ell}) \Big |^p_m  \Big ) \, .
$$
Applying Burkholder's inequality for multivariate martingales and using the  stationarity, we infer that there exists a positive constant $c_p$ depending only on $p$ such that, for any positive integer $n$,
$$
\Big \Vert \max_{r \leq k \leq n}  \Big | \sum_{\ell = r+1}^{k+r}  P_{\ell-j}(X_{\ell}) \Big |_m  \Big \Vert_p^{p'} \leq   c_p
 \sum_{\ell = r+1}^{n+r} \Vert  P_{\ell-j}(X_{\ell}) \Vert_{p,m}^{p'}  \leq  c_p \,
n \, \Vert  P_{0}(X_{j}) \Vert_{p,m}^{p'} \, .
$$
So, overall
\beq \label{decRn9}
\Big \Vert  \max_{r \leq k \leq n}  \Big | \sum_{\ell = r+1}^{k+r} \E_{\ell-r-1} (X_{\ell})  \Big |_m \Big \Vert^{p}_{p} \leq  ( c_p n)^{p/p'}
\sum_{j=r+1}^{\infty
} \alpha_j^{1-p} u_j^p = ( c_p n)^{p/p'}
\Big ( \sum_{j=r+1}^{\infty
} \Vert P_0(X_{j}) \Vert_p  \Big )^{p}\, .
\eeq
We handle now the second term in the right-hand side of \eqref{decRn8}. We first write that
\beq \label{decRn10}
S_k - \E_k (S_{k+r}) = S_k - S_{k-r} - \E_k ( S_{k+r} -S_{k-r}) + S_{k-r} -  \E_k (S_{k-r})\, .
\eeq
Let
$$
Y_r = \sum_{i=-(r-1)}^{0} X_i - \sum_{i=-(r-1)}^{r} \E_0(X_i) \, .
$$
With this notation, $$ S_k - S_{k-r} - \E_k ( S_{k+r} -S_{k-r}) = Y_r \circ T^k \, . $$ Hence, for any positive real $A$,
\begin{align*}
\Big \Vert \max_{r \leq k \leq n}  \big | S_k - S_{k-r} - \E_k ( S_{k+r} -S_{k-r})  \big |_m \Big \Vert^p_p & \leq 2^p A^p + 2^p
\big \Vert \max_{r \leq k \leq n}   | Y_r {\bf 1}_{|Y_r|_m >A} \circ T^k   |_m \big \Vert^p_p \\
& \leq 2^p A^p + 2^p \, n \, \Vert Y_r {\bf 1}_{|Y_r|_m >A}  \Vert^p_{p,m} \, .
\end{align*}
Since $ \Vert Y_r  \Vert_{p,m}  \leq K_r$ where $K_r$ is a constant depending on $r$, we get that
\beq \label{decRn11}
\lim_{A \rightarrow \infty} \limsup_{n \rightarrow \infty} \frac 1n\Big \Vert \max_{r \leq k \leq n}  \big | S_k - S_{k-r} - \E_k ( S_{k+r} -S_{k-r})  \big |_m \Big \Vert^{p'}_p  =0 \, .
\eeq
We deal now with the term $\Vert \max_{r \leq k \leq n} |S_{k-r} -  \E_k (S_{k-r}) |_m \Vert_p$. Since $\E_{\infty} (X_0) = X_0$ almost surely, we have that, almost surely
\begin{align*}
S_{k-r} -  \E_k (S_{k-r})  = \sum_{\ell=1}^{k-r} \sum_{j = - \infty}^{\ell-k-1} P_{\ell - j} (X_{\ell} )  = \sum_{j = - \infty}^{-k}  \sum_{\ell=1}^{k-r} P_{\ell - j} (X_{\ell} ) + \sum_{j = -k +1}^{-r-1}  \sum_{\ell= k+j+1}^{k-r} P_{\ell - j} (X_{\ell} ) \, .
\end{align*}
Therefore
\begin{multline} \label{decRn12}
\max_{r \leq k \leq n}  |S_{k-r} -  \E_k (S_{k-r}) |_m
\leq  \sum_{j = - \infty}^{-r}   \max_{r \leq k \leq n}  \Big | \sum_{\ell=1}^{k-r}  P_{\ell - j} (X_{\ell} )  \Big |_m + \sum_{j = -n +1}^{-r-1}   \max_{r \leq k \leq n}  \Big |  \sum_{\ell= k+j+1}^{k-r} P_{\ell - j} (X_{\ell} )  \Big |_m  \\
 \leq  \sum_{j = - \infty}^{-r}   \max_{r \leq k \leq n}  \Big | \sum_{\ell=1}^{k-r}  P_{\ell - j} (X_{\ell} )  \Big |_m + \sum_{j = -n +1}^{-r-1}   \max_{r \leq k \leq n}  \Big |  \sum_{\ell= r-n}^{k-r} P_{\ell - j} (X_{\ell} )  \Big |_m \\ + \sum_{j = -n +1}^{-r-1}   \max_{r \leq k \leq n}  \Big |  \sum_{\ell= r-n}^{k+j} P_{\ell - j} (X_{\ell} )  \Big |_m  \, .
\end{multline}
Let $u_{i} = \Vert P_0(X_{i}) \Vert_{p,m}$, $C_r = \sum_{ i=- \infty}^{-r}u_{i}$ and $\alpha_{i } = C_r^{-1} u_{i}$. As before, using the facts that for any $p \geq 1$, $x \mapsto x^p$ is convex and that $\alpha_{i}\geq
0$ with $\sum_{ i=- \infty}^{-r} \alpha_{i}=1$, we obtain that
$$
\Big \Vert  \sum_{j = - \infty}^{-r}   \max_{r \leq k \leq n}  \Big | \sum_{\ell=1}^{k-r}  P_{\ell - j} (X_{\ell} )  \Big |_m  \Big \Vert^{p}_{p}
\leq   \sum_{j = - \infty}^{-r}  \alpha_j^{1-p} \E \Big ( \max_{r \leq k \leq n}   \Big | \sum_{\ell=1}^{k-r}  P_{\ell - j} (X_{\ell} )  \Big |^p_m  \Big ) \, .
$$
Applying Burkholder's inequality for multivariate martingales
and using the stationarity, we infer that there exists a positive constant $c_p$ depending only on $p$ such that, for any positive integer $n$,
\beq \label{decRn13}
\Big \Vert  \sum_{j = - \infty}^{-r}   \max_{r \leq k \leq n}  \Big | \sum_{\ell=1}^{k-r}  P_{\ell - j} (X_{\ell} )  \Big |_m  \Big \Vert^{p}_{p}
\leq  (c_p n)^{p/p'} \Big ( \sum_{ i \geq r} \Vert P_0(X_{-i}) \Vert_p\Big )^p \, .
\eeq
With similar arguments, we derive that
\begin{multline}\label{decRn14}
\Big \Vert  \sum_{j = -n +1}^{-r-1}   \max_{r \leq k \leq n}   \Big |  \sum_{\ell= r-n}^{k-r} P_{\ell - j} (X_{\ell} )  \Big |_m \Big \Vert^{p}_{p}
 + \Big \Vert \sum_{j = -n +1}^{-r-1}   \max_{r \leq k \leq n}  \Big |  \sum_{\ell= r-n}^{k+j} P_{\ell - j} (X_{\ell} )  \Big |_m \Big \Vert^p_p
\\  \leq  2 (2 c_p n)^{p/p'} \Big ( \sum_{ i = r+1}^{n-1} \Vert P_0(X_{-i}) \Vert_p\Big )^p \, .
\end{multline}
Starting from \eqref{decRn12} and considering the upper bounds \eqref{decRn13} and \eqref{decRn14}, we get that
\beq \label{decRn15}
\Big \Vert  \max_{r \leq k \leq n}  |S_{k-r} -  \E_k (S_{k-r}) |_m  \Big \Vert_{p}
\leq  n^{1/p'} \sum_{ i \geq r} \Vert P_0(X_{-i}) \Vert_p \, .
\eeq
From the decomposition \eqref{decRn10}  together with \eqref{decRn11} and  \eqref{decRn15},
it follows that
\beq \label{decRn16}
\Big \Vert  \max_{r \leq k \leq n}  |S_{k} -  \E_k (S_{k+r}) |_m  \Big \Vert_{p}
\leq  n^{1/p'} \sum_{ i \geq r} \Vert P_0(X_{-i}) \Vert_p + o(n^{1/p'})\, .
\eeq
Starting from \eqref{decRn8} and considering \eqref{Burkh3bis}, \eqref{decRn9}, \eqref{decRn16} and the condition \eqref{Gor}, we derive that
\[
 \Big \Vert  \max_{r \leq k \leq n}   | R_k |_m
 \Big \Vert_{p} \ll n^{1/p'}  \sum_{|j| \geq r}
\Vert P_0(X_j)  \Vert_{p,m} + o(n^{1/p'})+r^{1/p'} \, ,
\]
(with the decomposition of   \eqref{Burkh3bis} in $\sum_{i=1}^r+\sum_{i=r+1}^n$)
which, combined with \eqref{app1} and Condition \eqref{Gor}, implies that

$$
\Big\Vert  \max_{1 \leq k \leq n}   | R_k |_m   \Big \Vert_{p}  \leq
\Big \Vert  \max_{1 \leq k \leq r}   | R_k |_m  \Big\Vert_p  +
\Big\Vert  \max_{r \leq k \leq n}   | R_k |_m    \Big\Vert_{p}
 \ll r^2 + n^{1/p'}  \sum_{|j| \geq r} \Vert P_0(X_j)  \Vert_{p,m} + o(n^{1/p'}) \, .
$$
Letting first $n$ tend to infinity and next  $r$ tend  to infinity, \eqref{app2} follows. $\square$

\medskip

Starting from Proposition \ref{approxrnp} one can prove the following theorem concerning the
weak and strong invariance principles for non-adapted sequences.

\begin{Theorem} \label{ASIPgene}
Let $X_0$ be a zero mean random variable in ${\mathbb L}^{2} ({\mathbb R}^m)$ and ${\cal F}_0 $ a $\sigma$-algebra satisfying ${\cal F}_0
\subseteq T^{-1 }({\cal F}_0)$. For any $i \in {\mathbb Z}$, let $X_i = X_0\circ T^i$ and ${\cal F}_i =T^{-i
}({\cal F}_0)$. Let $S_n =X_1 + \dots +X_n$. Assume that $T$ is ergodic, that $\E_{- \infty} (X_0) = 0$ almost surely, and that $\E_{\infty} (X_0) = X_0$
almost surely.
\begin{itemize}\item[1.] Assume that
\beq \label{condasipgene1}
\sum_{n \in {\mathbb Z}} \|P_0(X_n)\|_{2,m}<\infty \, .
\eeq
Then $n^{-1} {\rm Var} (S_n)$ converges to
\beq \label{defC}
C = \sum_{k \in {\mathbb Z}} {\mathrm{Cov}}\big (  X_0 , X_k \big ) \, .
\eeq
In addition the process $ \{ n^{-1/2} S_{[nt]}, t \in [0,1] \}$ converges in $D([0,1], {\mathbb R}^m )$ equipped with the uniform topology to a Wiener process $\{W(t), t\in [0,1]\}$ with variance matrix ${\rm Var} (W(1)) =C$.
\item[2.] Assume that
\beq \label{condasipgene}
\sum_{n\ge 3} \frac{\log n (\|P_0(X_n)\|_{2,m} + \|P_0(X_{-n})\|_{2,m} )}{(\log \log n)^{1/2}}<\infty \, .
\eeq
Then, enlarging the probability space if
necessary, there exists a sequence $(Z_i)_{i
\geq 1}$ of iid Gaussian random variables in ${\mathbb R}^m$ with zero mean and
variance matrix ${\rm Var}(Z_i)=C$ given by \eqref{defC}, such that \beq \label{resasipgene}
\sup_{1\leq k \leq n} \Big |\sum_{i=1}^k X_0 \circ T^i - \sum_{i=1}^k Z_i \Big |_m  = o \big (  ( n \log \log n)^{1/2}   \big )
\text{ almost surely, as $n\rightarrow \infty$}.
\eeq
\end{itemize}
\end{Theorem}

\begin{Remark}
The weak invariance principle (Item 1 of Theorem \ref{ASIPgene}) still holds if $T$ is not ergodic, but in that
case the limiting distribution is a mixture of Brownian motion (this has been
proved in \cite{DMV}
when $m=1$). This weak invariance
principle can be also extended to separable Hilbert spaces, with the appropriate
covariance operator. In the adapted case (i.e. $X_0$ is ${\mathcal F}_0$-measurable), the non ergodic Hilbert-valued version of
Item 1 has been proved in \cite{DM}.
\end{Remark}
\noindent{\bf Proof of Theorem \ref{ASIPgene}.} Let
$S_{m,1}=\{x \in {\mathbb R}^m: |x|_m = 1\}$. For a matrix $A$ from ${\mathbb R}^m$
to ${\mathbb R}^m$, let $|A|_m=\sup_{x \in S_{m,1}}|Ax|_m$.
By stationarity
 $n^{-1}{\rm Var} (S_n)=n^{-1}\sum_{|k|<n}(n-|k|){\mathrm{Cov}}(X_0,X_k)$. Hence $n^{-1}{\rm Cov} (S_n)$ converges to $C$ provided that
$ \sum_{k \in {\mathbb Z}} |{\rm Cov}(X_0, X_k)|_m < \infty$. Since $\E_{- \infty} (X_k) = 0$ almost surely and $\E_{\infty} (X_k) = X_k$ almost surely,
it follows that
$X_k=\sum_{i \in {\mathbb Z}}P_i(X_k)$ almost surely. Moreover
${\rm Cov}(P_i(X_0), P_j(X_k))=0$ for $i\neq j$. Hence,
$$
{\rm Cov}(X_0, X_k)=  \sum_{i \in {\mathbb Z}}
{\rm Cov}(P_i(X_0), P_i(X_k)) \, ,
$$
and consequently
$
|{\rm Cov}(X_0, X_k)|_m \leq \sum_{i \in {\mathbb Z}}
\|P_i(X_0)\|_{2,m}\|P_i(X_k)\|_{2,m} \, .
$
By \eqref{condasipgene1} it follows that
$$
  \sum_{k \in {\mathbb Z}} |{\rm Cov}(X_0, X_k)|_m \leq
  \sum_{k \in {\mathbb Z}}\sum_{i \in {\mathbb Z}} \|P_i(X_0)\|_{2,m}\|P_i(X_k)\|_{2,m}
  = \Big( \sum_{i \in {\mathbb Z}} \|P_0(X_i)\|_{2,m}\Big)^2 < \infty \, ,
$$ which proves the convergence of $n^{-1}{\rm Var} (S_n)$ to $C$.

Let now $d_0 := \sum_{j \in {\mathbb Z}} P_0 (X_j)$. Since \eqref{condasipgene1} is assumed, $d_0$ belongs to ${\mathbb L}^{2} ({\mathbb R}^m)$. In addition  $\E(d_0 |{\mathcal F}_{-1})=0$ almost surely.
Let $d_i := d_0 \circ T^i$ for all $i \in {\mathbb Z}$. Then $(d_i)_{i \in {\mathbb Z}}$ is a stationary ergodic sequence of martingale differences in ${\mathbb L}^{2} ({\mathbb R}^m)$. Let
$$ M_n := \sum_{i=1}^n d_i \, \text{ and } \, R_n := S_n -M_n \, .$$
Using \eqref{condasipgene1}, it follows from \eqref{app2} of Lemma \ref{approxrnp} that
\beq \label{p1item1}
\Big \|\max_{1 \leq k \leq n}| R_k|_m  \Big \|^{2}_{2} = o(n)  \, .
\eeq Since $n^{-1} {\rm Var} (S_n)$ converges to
$C$, it follows that ${\rm Var} (d_0)=C$. Therefore, Item 1 of Theorem \ref{ASIPgene} follows from the weak
invariance principle for partial sums  of  stationary  multivariate
martingale differences in ${\mathbb L}^{2} ({\mathbb R}^m)$ (see \cite{DM} for the non ergodic Hilbert-valued version) together with the maximal martingale
approximation given in \eqref{p1item1}.

\medskip

We turn now to the proof of Item 2. According to Theorem 3.1 in \cite{Be}
(that is the generalization of the Strassen's invariance principle \cite{ST} for
real martingales with ergodic increments to the multivariate case),
enlarging the probability space if
necessary, there exists a sequence $(Z_i)_{i
\geq 1}$ of iid Gaussian random variables in ${\mathbb R}^m$ with zero mean and
covariance ${\mathrm{Var}}(Z_1)=C $ such that $$
\sup_{1\leq k \leq n} \Big |\sum_{i=1}^k d_0\circ T^i - \sum_{i=1}^k Z_i \Big |_m  = o \big ( ( n \log \log n)^{1/2} \big ) \
\text{ almost surely, as $n\rightarrow \infty$}.
$$
Therefore the strong approximation result \eqref{resasipgene} will follow if we can show that
\beq \label{convrnps}
|R_n |_m  = o \big ( ( n \log \log n)^{1/2}  \big ) \
\text{ almost surely, as $n\rightarrow \infty$}.
\eeq
Since $R_n = \sum_{i=1}^n (f - d_0)\circ T^i$, \eqref{convrnps} will follow by Theorem 4.7 in \cite{CM} if we can prove that
$$
\sum_{n >3} \frac{\Vert R_n \Vert_2}{n^{3/2} ( \log \log  n)^{1/2} } < \infty \, .
$$
Using \eqref{app1} of Lemma \ref{approxrnp}, this last convergence will hold provided that
\beq \label{convrnps2}
\sum_{n >3} \frac{\Big ( \sum_{k= 1}^{n} \Big ( \sum_{|\ell| \ge k}\| P_{\ell}(X_0)\|_{2,m} \Big )^{2} \Big )^{1/2}}{n^{3/2} ( \log \log  n)^{1/2}} < \infty \, .
\eeq
Notice that
\begin{multline*}
 \sum_{n >3}  \frac{\Big ( \sum_{k= 1}^{n} \Big ( \sum_{|\ell| \ge k}  \| P_{\ell}(X_0)\|_{2,m} \Big )^{2} \Big )^{1/2}}{n^{3/2} ( \log \log  n)^{1/2}} \\
  \ll
\sum_{k >2} 2^{-k/2} ( \log k )^{-1/2} \Big ( \sum_{j=1}^{2^k} \Big ( \sum_{\ell \geq j}  ( \|P_0(X_\ell)\|_{2,m} +  \|P_0(X_{-\ell})\|_{2,m} ) \Big )^2 \Big ) ^{1/2} \\
 \ll \sum_{k \geq 0}2^{-k/2} ( \log k )^{-1/2} \Big (  \sum_{j=0}^{k} 2^j\Big ( \sum_{\ell
\geq 2^j} ( \|P_0(X_\ell)\|_{2,m} +  \|P_0(X_{-\ell})\|_{2,m} ) \Big )^2 \Big ) ^{1/2}\, .
\end{multline*}
Now, using the subadditivity of $x\mapsto x^{1/2}$, it follows that \eqref{convrnps2} will be satisfied as soon as
$$
\sum_{k >2}2^{-k/2} ( \log k )^{-1/2}  \sum_{j=0}^{k} 2^{j/2} \sum_{\ell \geq 2^j} ( \|P_0(X_\ell)\|_{2,m} +  \|P_0(X_{-\ell})\|_{2,m} )  < \infty \, ,
$$
which holds as soon as  \eqref{condasipgene} does (changing the order of summation
in $\sum_\ell\sum_j\sum_k$). This ends the proof of Item 2 of Theorem
\ref{ASIPgene}. $\square$

\medskip

For the sake of applications, we now give  sufficient conditions
for \eqref{condasipgene1} and \eqref{condasipgene} to hold.
\begin{Remark} \label{use}
The condition \eqref{condasipgene1} is satisfied if we assume that
\begin{equation}  \label{usecond1}
\sum_{n\geq 1 }\frac{1}{n^{1/2}}\Vert \E_0 (X_n)\Vert _{2,m}<\infty \text{ and }\sum_{n\geq 1 }\frac{1}{n^{1/2}}%
\Vert X_{-n}-\E_0 (X_{-n}) \Vert _{2,m}<\infty \, ,
\end{equation}
and the condition \eqref{condasipgene} holds if we assume that
\begin{equation}  \label{usecond2}
\sum_{n\geq 3 }\frac{\log n}{n^{1/2}(\log \log n)^{1/2}}\Vert \E_0 (X_n)\Vert _{2,m}<\infty \text{ and }\sum_{n\geq 3 }\frac{\log n}{n^{1/2}(\log \log n)^{1/2}}%
\Vert X_{-n}-\E_0 (X_{-n}) \Vert _{2,m}<\infty \, .
\end{equation}
\end{Remark}
The proof of the remark above is omitted since it
uses exactly the arguments developed to prove Remarks 3.3 and 3.6 in \cite{DeMePe} (see Section 5.5 of \cite{DeMePe}). Notice that the conditions \eqref{usecond1} or \eqref{usecond2} imply clearly that $ \E_{-\infty}(X_0) = 0$ almost surely and that $\E_{\infty} (X_0) = X_0$ almost surely.

\end{document}